\documentclass[preprint,10pt]{elsarticle}
\usepackage{color}
\usepackage{geometry}
\usepackage{amssymb}
\usepackage{booktabs}
\usepackage{rotating}
\usepackage{amsmath}
\usepackage{amsfonts} 

\newtheorem{example}{Example}[section]
\newtheorem{remark}{Remark}

\newtheorem{theorem}{Theorem}[section]

\biboptions{numbers,sort&compress}
\journal{arXiv.org}

\begin{document}

\begin{frontmatter}
\title{A new $Z$-eigenvalue inclusion theorem for tensors\tnoteref{label1}}
\tnotetext[label1]{This work is supported by the National Natural Science Foundations of China (Grant Nos.11361074,11501141),
Foundation of Guizhou Science and Technology Department (Grant No.[2015]2073)
and Natural Science Programs of Education Department of Guizhou Province (Grant No.[2016]066).}


\author{Jianxing Zhao\corref{cor1}}
\ead{zjx810204@163.com}
\address{College of Data Science and Information Engineering, Guizhou Minzu University,
Guiyang 550025, P.R.China}
\cortext[cor1]{Corresponding author.}
\begin{abstract}
A new $Z$-eigenvalue inclusion theorem for tensors is given and proved to be tighter than those in
[G. Wang, G.L. Zhou, L. Caccetta, $Z$-eigenvalue inclusion theorems for tensors, Discrete and Continuous Dynamical Systems Series B,
22(1) (2017) 187--198]. Based on this set, a sharper upper bound for the $Z$-spectral radius of weakly symmetric nonnegative tensors is obtained. Finally, numerical examples are given to show the effectiveness of the proposed bound.
\end{abstract}

\end{frontmatter}


\section{Introduction}
For a positive integer $n$,~$n\geq 2$, $N$ denotes the set $\{1,2,\cdots,n\}$.
$\mathbb{C}$ ($\mathbb{R}$) denotes the set of all complex (real) numbers.
We call $\mathcal{A}=(a_{i_1i_2\cdots i_{m}})$ a real tensor of order $m$ dimension $n$, denoted by $\mathbb{R}^{[m,n]}$, if
$$a_{i_{1}i_2\cdots i_{m}}\in{\mathbb{R}},$$ where $i_{j}\in{N}$ for $j=1,2,\cdots,m$.
$\mathcal{A}$ is called nonnegative if $a_{i_{1}i_2\cdots i_{m}}\geq 0.$
$\mathcal{A}=(a_{i_{1}\cdots i_{m}})\in \mathbb{R}^{[m,n]}$ is called symmetric \cite{qi2005} if
\begin{eqnarray}
a_{i_{1}\cdots i_{m}}=a_{\pi(i_{1}\cdots i_{m})},\ \forall \pi\in\Pi_{m},\nonumber
\end{eqnarray}
where $\Pi_{m}$ is the permutation group of $m$ indices.
$\mathcal{A}=(a_{i_1\cdots i_{m}})\in \mathbb{R}^{[m,n]}$ is called weakly symmetric \cite{changkc} if the associated homogeneous polynomial
$$\mathcal{A}x^m=\sum\limits_{i_1,\cdots,i_m\in N}a_{i_1\cdots i_{m}}x_{i_1}\cdots x_{i_m}$$
satisfies $\nabla \mathcal{A}x^m=m\mathcal{A}x^{m-1}$.
It is shown in \cite{changkc} that a symmetric tensor is necessarily weakly symmetric, but the converse
is not true in general.

Given a tensor $\mathcal {A}=(a_{i_1\cdots i_m})\in \mathbb{R}^{[m,n]}$, if there are $\lambda\in \mathbb{C}$ and
$x=(x_1,x_{2}\cdots,x_n)^T\in \mathbb{C}\backslash\{0\}$
such that
\[\mathcal {A}x^{m-1}=\lambda x~\textmd{and}~ x^Tx=1,\]
then $\lambda$ is called an $E$-eigenvalue of $\mathcal {A}$ and $x$ an
$E$-eigenvector of $\mathcal {A}$ associated with $\lambda$, where
$\mathcal {A}x^{m-1}$ is an $n$ dimension vector whose $i$th component is
\[(\mathcal
{A}x^{m-1})_i=\sum\limits_{i_2,\cdots,i_m\in N} a_{ii_2\cdots
i_m}x_{i_2}\cdots x_{i_m}.\]
If $\lambda$ and $x$ are all real, then $\lambda$ is called a $Z$-eigenvalue of $\mathcal
{A}$ and $x$ a $Z$-eigenvector of $\mathcal {A}$ associated with
$\lambda$; for details, see \cite{qi2005,lim}.

We define the $Z$-spectrum of $\mathcal{A}$,
denoted $\sigma(\mathcal{A})$ to be the set of all $Z$-eigenvalues of $\mathcal{A}$. Assume $\sigma(\mathcal{A})\neq 0,$ then the $Z$-spectral radius \cite{changkc} of $\mathcal{A}$, denoted $\varrho(\mathcal{A})$, is defined as
$$\varrho(\mathcal{A}):=\sup\{|\lambda|:\lambda\in \sigma(\mathcal{A})\}.$$

Recently, much literature has focused on locating all $Z$-eigenvalues of tensors and
bounding the $Z$-spectral radius of nonnegative tensors in \cite{wg,sys,liwen,hejunjcaa,hejunspringerplus,hjaml,lql}.
It is well known that one can use eigenvalue inclusion sets to obtain the lower and upper bounds of the spectral radius of nonnegative tensors; for details, see \cite{wg,lcq-lyt,lcq-kx,lcq-zjj,lcq-cz}. Therefore, the main aim of this paper is to give a tighter $Z$-eigenvalue inclusion set for tensors, and use it to obtain a sharper upper bound for the $Z$-spectral radius of weakly symmetric nonnegative tensors.

In 2017, Wang \emph{et al}. \cite{wg} established the following Ger$\breve{s}$gorin-type  $Z$-eigenvalue inclusion theorem for tensors.
\begin{theorem}\emph{\cite[Theorem 3.1]{wg}}\label{wg-th1}
Let $\mathcal{A}=(a_{i_{1}\cdots i_{m}})\in{\mathbb{R}}^{[m,n]}$. Then
\begin{eqnarray*}
\sigma(\mathcal{A})\subseteq \mathcal{K}(\mathcal{A})=\bigcup\limits_{i\in{N}}\mathcal{K}_{i}(\mathcal{A}),
\end{eqnarray*}
where 
\begin{eqnarray*}
\mathcal{K}_{i}(\mathcal{A})=\{z\in{\mathbb{C}}:|z|\leq R_{i}(\mathcal{A})\},
~R_{i}(\mathcal{A})=\sum\limits_{i_2,\cdots, i_m\in N}|a_{ii_{2}\cdots i_{m}}|.
\end{eqnarray*}
\end{theorem}

To get tighter $Z$-eigenvalue inclusion sets than $\mathcal{K}(\mathcal{A})$,
Wang \emph{et al}. \cite{wg} also gave a Brauer-type $Z$-eigenvalue inclusion theorem for tensors.

\begin{theorem}\emph{\cite[Theorem 3.3]{wg}}\label{wg-th2}
Let $\mathcal{A}=(a_{i_{1}\cdots i_{m}})\in{\mathbb{R}}^{[m,n]}$. Then
\begin{eqnarray*}
\sigma(\mathcal{A})\subseteq\mathcal{M}(\mathcal{A})=\bigcup\limits_{i,j\in N,i\neq j}
\left(\mathcal{M}_{i,j}(\mathcal{A})\bigcup\mathcal{H}_{i,j}(\mathcal{A})\right),
\end{eqnarray*}
where
\[
\mathcal{M}_{i,j}(\mathcal{A})=\left\{z\in{\mathbb{C}}:
\big[|z|-(R_i(\mathcal{A})-|a_{ij\cdots j}|\big](|z|-P_j^i(\mathcal{A}))\leq |a_{ij\cdots j}|(R_j(\mathcal{A})-P_j^i(\mathcal{A}))\right\},
\]
\[
\mathcal{H}_{i,j}(\mathcal{A})=\left\{z\in{\mathbb{C}}: |z|<R_i(\mathcal{A})-|a_{ij\cdots j}|, |z|< P_j^i(\mathcal{A})\right\},
\]
and
\[P_j^i(\mathcal{A})=\sum\limits_{i_2,\cdots, i_m\in N,\atop i\notin\{i_2,\cdots,i_m\}}|a_{ji_{2}\cdots i_{m}}|.\]
\end{theorem}

In this paper, we continue this research on the $Z$-eigenvalue localization problem for tensors and its applications.
We give a new $Z$-eigenvalue inclusion set for tensors and prove that the new set is tighter than those in
Theorem \ref{wg-th1} and Theorem \ref{wg-th2}.
As an application of this set, we obtain a new upper bound for the $Z$-spectral radius of weakly symmetric nonnegative tensors, which is sharper than existing bounds in some cases.

\section{A new $Z$-eigenvalue inclusion theorem}\label{sec2}
In this section, we give a new $Z$-eigenvalue inclusion theorem for tensors, and establish the comparison between this set with those in Theorem \ref{wg-th1} and Theorem \ref{wg-th2}.

\begin{theorem}\label{th1}
Let $\mathcal{A}=(a_{i_{1}\cdots i_{m}})\in{\mathbb{R}}^{[m,n]}$. Then
\begin{eqnarray*}
\sigma(\mathcal{A})\subseteq \Omega(\mathcal{A})=
\bigcup\limits_{i,j\in N, j\neq i}\Bigg(\hat{\Omega}_{i,j}(\mathcal{A})\bigcup\Big(\tilde{\Omega}_{i,j}(\mathcal{A})\bigcap \mathcal{K}_i(\mathcal{A})\Big)\Bigg),
\end{eqnarray*}
where
\begin{eqnarray*}
\hat{\Omega}_{i,j}(\mathcal{A})=\left\{z\in \mathbb{C}:
|z|< P_i^j(\mathcal{A}), |z|< P_j^i(\mathcal{A})\right\}
\end{eqnarray*}
and
\begin{eqnarray*}
\tilde{\Omega}_{i,j}(\mathcal{A})=\left\{z\in \mathbb{C}:
\big(|z|-P_i^j(\mathcal{A})\big)\big(|z|-P_j^i(\mathcal{A})\big)
\leq \big(R_i(\mathcal{A})-P_i^j(\mathcal{A})\big)\big(R_j(\mathcal{A})-P_j^i(\mathcal{A})\big)
\right\}.
\end{eqnarray*}
\end{theorem}
\noindent {\bf Proof.} 
Let $\lambda$ be a $Z$-eigenvalue of $\mathcal{A}$ with corresponding $Z$-eigenvector $x=(x_{1},\cdots,x_{n})^{T}\in{\mathbb{C}}^{n}\backslash\{0\}$, i.e.,
\begin{eqnarray}\label{th1-equ1}
\mathcal{A}x^{m-1}=\lambda x,~\textmd{and}~||x||_2=1.
\end{eqnarray}
Let $|x_t|\geq |x_s|\geq\max\limits_{i \in N,i\neq t,s}|x_i|$. Obviously, $0<|x_t|^{m-1}\leq |x_t|\leq 1.$
From (\ref{th1-equ1}), we have
\begin{eqnarray*}
\lambda x_t=\sum\limits_{i_2,\cdots, i_m\in N,\atop s\in \{i_2,\cdots i_m\}}a_{ti_{2}\cdots i_{m}}x_{i_{2}}\cdots x_{i_{m}}
+\sum\limits_{i_2,\cdots, i_m\in N,\atop s\notin \{i_2,\cdots i_m\}}a_{ti_{2}\cdots i_{m}}x_{i_{2}}\cdots x_{i_{m}}.
\end{eqnarray*}
Taking modulus in the above equation and using the triangle inequality gives
\begin{eqnarray*}
|\lambda||x_t|&\leq& \sum\limits_{i_2,\cdots, i_m\in N,\atop s\in \{i_2,\cdots i_m\}}|a_{ti_{2}\cdots i_{m}}||x_{i_{2}}|\cdots |x_{i_{m}}|
+\sum\limits_{i_2,\cdots, i_m\in N,\atop s\notin \{i_2,\cdots i_m\}}|a_{ti_{2}\cdots i_{m}}||x_{i_{2}}|\cdots |x_{i_{m}}|\\
&\leq& \sum\limits_{i_2,\cdots, i_m\in N,\atop s\in \{i_2,\cdots i_m\}}|a_{ti_{2}\cdots i_{m}}||x_s|
+\sum\limits_{i_2,\cdots, i_m\in N,\atop s\notin \{i_2,\cdots i_m\}}|a_{ti_{2}\cdots i_{m}}||x_t|\\
&=&(R_t(\mathcal{A})-P_t^s(\mathcal{A}))|x_s|+P_t^s(\mathcal{A})|x_t|,
\end{eqnarray*}
i.e.,
\begin{eqnarray}\label{th1-equ2}
\big(|\lambda|-P_t^s(\mathcal{A})\big)|x_t|\leq (R_t(\mathcal{A})-P_t^s(\mathcal{A}))|x_s|.
\end{eqnarray}
If $|x_s|=0$, then $|\lambda|-P_t^s(\mathcal{A})\leq 0$ as $|x_t|>0$.
When $|\lambda|\geq P_s^t(\mathcal{A})$, we have
\begin{eqnarray*}
\big(|\lambda|-P_t^s(\mathcal{A})\big)\big(|\lambda|-P_s^t(\mathcal{A})\big)
\leq 0\leq (R_t(\mathcal{A})-P_t^s(\mathcal{A}))\big(R_s(\mathcal{A})-P_s^t(\mathcal{A})\big),
\end{eqnarray*}
which implies $\lambda\in \tilde{\Omega}_{t,s}(\mathcal{A})\subseteq \Omega(\mathcal{A})$.
When $|\lambda|< P_s^t(\mathcal{A})$, we have
$\lambda\in \hat{\Omega}_{t,s}(\mathcal{A})\subseteq \Omega(\mathcal{A})$.

Otherwise, $|x_s|>0$. By (\ref{th1-equ1}), we can get
\begin{eqnarray*}
|\lambda| |x_s|&\leq&\sum\limits_{i_2,\cdots, i_m\in N,\atop t\in \{i_2,\cdots i_m\}}|a_{si_{2}\cdots i_{m}}||x_{i_{2}}|\cdots |x_{i_{m}}|
+\sum\limits_{i_2,\cdots, i_m\in N,\atop t\notin \{i_2,\cdots i_m\}}|a_{si_{2}\cdots i_{m}}||x_{i_{2}}|\cdots |x_{i_{m}}|\\
&\leq&\sum\limits_{i_2,\cdots, i_m\in N,\atop t\in \{i_2,\cdots i_m\}}|a_{si_{2}\cdots i_{m}}||x_t|^{m-1}
+\sum\limits_{i_2,\cdots, i_m\in N,\atop t\notin (i_2,\cdots i_m)}|a_{si_{2}\cdots i_{m}}||x_s|^{m-1},\\
&\leq&\sum\limits_{i_2,\cdots, i_m\in N,\atop t\in \{i_2,\cdots i_m\}}|a_{si_{2}\cdots i_{m}}||x_t|
+\sum\limits_{i_2,\cdots, i_m\in N,\atop t\notin (i_2,\cdots i_m)}|a_{si_{2}\cdots i_{m}}||x_s|,
\end{eqnarray*}
i.e.,
\begin{eqnarray}\label{th1-equ3}
\big(|\lambda|-P_s^t(\mathcal{A})\big) |x_s|\leq (R_s(\mathcal{A})-P_s^t(\mathcal{A}))|x_t|.
\end{eqnarray}
By (\ref{th1-equ2}), it is not difficult to see $|\lambda|\leq R_t(\mathcal{A})$, that is, $\lambda\in \mathcal{K}_t(\mathcal{A})$.
When $|\lambda|\geq P_t^s(\mathcal{A})$ or $|\lambda|\geq P_s^t(\mathcal{A})$ holds,
multiplying (\ref{th1-equ2}) with (\ref{th1-equ3}) and noting that $|x_t||x_s|>0$, we have
\begin{eqnarray*}
\big(|\lambda|-P_t^s(\mathcal{A})\big)\big(|\lambda|-P_s^t(\mathcal{A})\big)
\leq \big(R_t(\mathcal{A})-P_t^s(\mathcal{A})\big)\big(R_s(\mathcal{A})-P_s^t(\mathcal{A})\big),
\end{eqnarray*}
which implies
$\lambda\in \big(\tilde{\Omega}_{t,s}(\mathcal{A})\bigcap \mathcal{K}_t(\mathcal{A})\big)\subseteq \Omega(\mathcal{A})$.

And when $|\lambda|< P_t^s(\mathcal{A})$ and $|\lambda|< P_s^t(\mathcal{A})$ hold,
we have $\lambda\in \hat{\Omega}_{t,s}(\mathcal{A})\subseteq \Omega(\mathcal{A})$.
Hence, the conclusion $\sigma(\mathcal{A})\subseteq \Omega(\mathcal{A})$ follows immediately from what we have proved.
\hfill$\blacksquare$ 

Next, a comparison theorem is given for Theorem \ref{wg-th1}, Theorem \ref{wg-th2} and Theorem \ref{th1}.
\begin{theorem}\label{th2}
Let $\mathcal{A}=(a_{i_{1}\cdots i_{m}})\in{\mathbb{R}}^{[m,n]}$. Then
\begin{eqnarray*}
\Omega(\mathcal{A})\subseteq \mathcal{M}(\mathcal{A})\subseteq \mathcal{K}(\mathcal{A}).
\end{eqnarray*}
\end{theorem}
\noindent {\bf Proof.}
By Corollary 3.2 in \cite{wg}, $\mathcal{M}(\mathcal{A})\subseteq \mathcal{K}(\mathcal{A})$ holds.
Hence, we only prove $\Omega(\mathcal{A})\subseteq \mathcal{M}(\mathcal{A})$.
Let $z\in\Omega(\mathcal{A})$.
Then there are $t,s\in N$ and $t\neq s$ such that
$z\in \hat{\Omega}_{t,s}(\mathcal{A})$ or $z\in \big(\tilde{\Omega}_{t,s}(\mathcal{A})\bigcap \mathcal{K}_t(\mathcal{A})\big)$.
We divide the proof into two parts.

Case I: If $z\in \hat{\Omega}_{t,s}(\mathcal{A}),$ that is,
$|z|< P_t^s(\mathcal{A})$ and $|z|< P_s^t(\mathcal{A})$.
Then, it is easily to see that
$$|z|< P_t^s(\mathcal{A})\leq R_t(\mathcal{A})-|a_{ts\cdots s}|,$$
which implies that
$z\in \mathcal{H}_{t,s}(\mathcal{A})\subseteq \mathcal{M}(\mathcal{A})$,
consequently, $\Omega(\mathcal{A})\subseteq \mathcal{M}(\mathcal{A})$.

Case II: If $z\notin \hat{\Omega}_{t,s}(\mathcal{A}),$ that is,
\begin{eqnarray}\label{th2-equ1}
|z|\geq P_s^t(\mathcal{A})
\end{eqnarray}
or
\begin{eqnarray}\label{th2-equ2}
|z|\geq P_t^s(\mathcal{A}),
\end{eqnarray}
then
$z\in \Big(\tilde{\Omega}_{t,s}(\mathcal{A})\bigcap \mathcal{K}_t(\mathcal{A})\Big)$, i.e.,
\begin{eqnarray}\label{th2-equ3}
|z|\leq R_t(\mathcal{A})
\end{eqnarray}
and
\begin{eqnarray}\label{th2-equ4}
\big(|z|-P_t^s(\mathcal{A})\big)\big(|z|-P_s^t(\mathcal{A})\big)
\leq \big(R_t(\mathcal{A})-P_t^s(\mathcal{A})\big)\big(R_s(\mathcal{A})-P_s^t(\mathcal{A})\big).
\end{eqnarray}

(i) Assume $\big(R_t(\mathcal{A})-P_t^s(\mathcal{A})\big)\big(R_s(\mathcal{A})-P_s^t(\mathcal{A})\big)=0$.
When (\ref{th2-equ1}) holds, we have
\begin{eqnarray*}
\big[|z|-(R_t(\mathcal{A})-|a_{ts\cdots s}|)\big]\big(|z|-P_s^t(\mathcal{A})\big)&\leq&\big(|z|-P_t^s(\mathcal{A})\big)\big(|z|-P_s^t(\mathcal{A})\big)\\
&\leq& \big(R_t(\mathcal{A})-P_t^s(\mathcal{A})\big)\big(R_s(\mathcal{A})-P_s^t(\mathcal{A})\big)\\
&=&0\leq |a_{ts\cdots s}|\big(R_s(\mathcal{A})-P_s^t(\mathcal{A})\big),
\end{eqnarray*}
which implies that $z\in \mathcal{M}_{t,s}(\mathcal{A})\subseteq \mathcal{M}(\mathcal{A}).$

On the other hand, when (\ref{th2-equ2}) holds and $|z|< P_s^t(\mathcal{A})$, we have
$z\in \mathcal{H}_{t,s}(\mathcal{A})\subseteq \mathcal{M}(\mathcal{A})$
if
$$P_t^s(\mathcal{A})\leq |z|< R_t(\mathcal{A})-|a_{ts\cdots s}|,$$
and $z\in \mathcal{M}_{t,s}(\mathcal{A})\subseteq \mathcal{M}(\mathcal{A})$ from
\begin{eqnarray*}
\big[|z|-(R_t(\mathcal{A})-|a_{ts\cdots s}|)\big]\big(|z|-P_s^t(\mathcal{A})\big)\leq 0
\leq |a_{ts\cdots s}|\big(R_s(\mathcal{A})-P_s^t(\mathcal{A})\big)
\end{eqnarray*}
if
$$R_t(\mathcal{A})-|a_{ts\cdots s}|\leq |z|\leq R_t(\mathcal{A}).$$

(ii) Assume $\big(R_t(\mathcal{A})-P_t^s(\mathcal{A})\big)\big(R_s(\mathcal{A})-P_s^t(\mathcal{A})\big)>0$.
Then dividing both sides by $\big(R_t(\mathcal{A})-P_t^s(\mathcal{A})\big)\big(R_s(\mathcal{A})-P_s^t(\mathcal{A})\big)$ in (\ref{th2-equ4}), we have
\begin{eqnarray}\label{th2-equ5}
\frac{|z|-P_t^s(\mathcal{A})}{R_t(\mathcal{A})-P_t^s(\mathcal{A})}
\frac{|z|-P_s^t(\mathcal{A})}{R_s(\mathcal{A})-P_s^t(\mathcal{A})}\leq 1.
\end{eqnarray}

Let $a=|z|,b=P_t^s(\mathcal{A}),c=R_t(\mathcal{A})-|a_{ts\cdots s}|-P_t^s(\mathcal{A})$ and $d=|a_{ts\cdots s}|$.
If $|a_{ts\cdots s}|>0$, by (\ref{th2-equ3}) and Lemma 2.2 in \cite{lcq-lyt}, we have
\begin{eqnarray}\label{th2-equ6}
\frac{|z|-(R_t(\mathcal{A})-|a_{ts\cdots s}|)}{|a_{ts\cdots s}|}=\frac{a-(b+c)}{d}
\leq \frac{a-b}{c+d}=\frac{|z|-P_t^s(\mathcal{A})}{R_t(\mathcal{A})-P_t^s(\mathcal{A})}.
\end{eqnarray}
When (\ref{th2-equ1}) holds, by (\ref{th2-equ5}) and (\ref{th2-equ6}), we have
\begin{eqnarray*}
\frac{|z|-(R_t(\mathcal{A})-|a_{ts\cdots s}|)}{|a_{ts\cdots s}|}\frac{|z|-P_s^t(\mathcal{A})}{R_s(\mathcal{A})-P_s^t(\mathcal{A})}
\leq \frac{|z|-P_t^s(\mathcal{A})}{R_t(\mathcal{A})-P_t^s(\mathcal{A})}\frac{|z|-P_s^t(\mathcal{A})}{R_s(\mathcal{A})-P_s^t(\mathcal{A})}
\leq 1,
\end{eqnarray*}
equivalently,
\begin{eqnarray*}
\big[|z|-(R_t(\mathcal{A})-|a_{ts\cdots s}|)\big]\big(|z|-P_s^t(\mathcal{A})\big)
\leq |a_{ts\cdots s}|\big(R_s(\mathcal{A})-P_s^t(\mathcal{A})\big),
\end{eqnarray*}
which implies that
$z\in \mathcal{M}_{t,s}(\mathcal{A})\subseteq \mathcal{M}(\mathcal{A})$.
On the other hand, when (\ref{th2-equ2}) holds and $|z|< P_s^t(\mathcal{A})$, we have
$z\in \mathcal{H}_{t,s}(\mathcal{A})\subseteq \mathcal{M}(\mathcal{A})$
if
$$P_t^s(\mathcal{A})\leq |z|< R_t(\mathcal{A})-|a_{ts\cdots s}|,$$
and $z\in \mathcal{M}_{t,s}(\mathcal{A})\subseteq \mathcal{M}(\mathcal{A})$ from
\begin{eqnarray*}
\big[|z|-(R_t(\mathcal{A})-|a_{ts\cdots s}|)\big]\big(|z|-P_s^t(\mathcal{A})\big)\leq 0
\leq |a_{ts\cdots s}|\big(R_s(\mathcal{A})-P_s^t(\mathcal{A})\big)
\end{eqnarray*}
if $R_t(\mathcal{A})-|a_{ts\cdots s}|\leq |z|\leq R_t(\mathcal{A}).$

If $|a_{ts\cdots s}|=0$, by $|z|\leq R_t(\mathcal{A})$, we have
\begin{eqnarray}\label{th2-equ7}
|z|-(R_t(\mathcal{A})-|a_{ts\cdots s}|)\leq 0=|a_{ts\cdots s}|.
\end{eqnarray}
When (\ref{th2-equ1}) holds, by (\ref{th2-equ7}), we can obtain
\begin{eqnarray*}
\big[|z|-(R_t(\mathcal{A})-|a_{ts\cdots s}|)\big]\big(|z|-P_s^t(\mathcal{A})\big)\leq 0= |a_{ts\cdots s}|\big(R_s(\mathcal{A})-P_s^t(\mathcal{A})\big),
\end{eqnarray*}
which implies that $z\in \mathcal{M}_{t,s}(\mathcal{A})\subseteq \mathcal{M}(\mathcal{A}).$
On the other hand, when (\ref{th2-equ2}) holds and $|z|< P_s^t(\mathcal{A})$, we easily get
$z\in \mathcal{H}_{t,s}(\mathcal{A})\subseteq \mathcal{M}(\mathcal{A})$
if $$P_t^s(\mathcal{A})\leq |z|< R_t(\mathcal{A})-|a_{ts\cdots s}|,$$
and $z\in \mathcal{M}_{t,s}(\mathcal{A})\subseteq \mathcal{M}(\mathcal{A})$ from
\begin{eqnarray*}
\big[|z|-(R_t(\mathcal{A})-|a_{ts\cdots s}|)\big]\big(|z|-P_s^t(\mathcal{A})\big)\leq 0= |a_{ts\cdots s}|\big(R_s(\mathcal{A})-P_s^t(\mathcal{A})\big)
\end{eqnarray*}
if $$R_t(\mathcal{A})-|a_{ts\cdots s}|\leq |z|\leq R_t(\mathcal{A}).$$
The conclusion follows from Case I and Case II.
\hfill$\blacksquare$

\begin{remark}\emph{
Theorem \ref{th2} shows that the set $\Omega(\mathcal{A})$ in Theorem \ref{th1} is tighter than $\mathcal{K}(\mathcal{A})$ in Theorem \ref{wg-th1}
and $\mathcal{M}(\mathcal{A})$ in Theorem \ref{wg-th2}, that is,
$\Omega(\mathcal{A})$ can capture all $Z$-eigenvalues of $\mathcal{A}$ more precisely than $\mathcal{K}(\mathcal{A})$ and $\mathcal{M}(\mathcal{A})$.
}\end{remark}

\section{A new upper bound for the $Z$-spectral radius of weakly symmetric nonnegative tensors}\label{sec3}
As an application of the results in Section 2, a new upper bound for the $Z$-spectral radius of weakly symmetric nonnegative tensors is given.

\begin{theorem}\label{th3}
Let $\mathcal{A}=(a_{i_{1}\cdots i_{m}})\in{\mathbb{R}}^{[m,n]}$ be a weakly symmetric nonnegative tensor. Then
\begin{eqnarray*}
\varrho(\mathcal{A})\leq \Omega_{max}=\max\big\{\hat{\Omega}_{max},\tilde{\Omega}_{max}\big\},
\end{eqnarray*}
where
\begin{eqnarray*}
\hat{\Omega}_{max}=\max\limits_{i,j\in N, j\neq i}\min\{P_i^j(\mathcal{A}),P_j^i(\mathcal{A})\},\\
\tilde{\Omega}_{max}=\max\limits_{i,j\in N,j\neq i}\min\left\{R_i(\mathcal{A}),\Delta_{i,j}(\mathcal{A})\right\},
\end{eqnarray*}
and
\begin{eqnarray*}
\Delta_{i,j}(\mathcal{A})=\frac{1}{2}\left\{
P_i^j(\mathcal{A})+P_j^i(\mathcal{A})+
\sqrt{\big(P_i^j(\mathcal{A})-P_j^i(\mathcal{A})\big)^2+4\big(R_i(\mathcal{A})-P_i^j(\mathcal{A})\big)\big(R_j(\mathcal{A})-P_j^i(\mathcal{A})\big)}
\right\}.
\end{eqnarray*}
\end{theorem}
\noindent {\bf Proof.}
From Lemma 4.4 in \cite{wg}, we know that $\varrho(\mathcal{A})$ is the largest $Z$-eigenvalue of $\mathcal{A}$.
By Theorem \ref{th1}, we have
\begin{eqnarray*}
\varrho(\mathcal{A})\in \bigcup\limits_{i,j\in N, j\neq i}\Bigg(\hat{\Omega}_{i,j}(\mathcal{A})\bigcup\Big(\tilde{\Omega}_{i,j}(\mathcal{A})\bigcap \mathcal{K}_i(\mathcal{A})\Big)\Bigg),
\end{eqnarray*}
that is, there are $t,s\in N, t\neq s$ such that $\varrho(\mathcal{A})\in\hat{\Omega}_{t,s}(\mathcal{A})$ or
$\varrho(\mathcal{A})\in \Big(\tilde{\Omega}_{t,s}(\mathcal{A})\bigcap \mathcal{K}_t(\mathcal{A})\Big)$.

If $\varrho(\mathcal{A})\in \hat{\Omega}_{t,s}(\mathcal{A}),$
i.e.,
$\varrho(\mathcal{A})< P_t^s(\mathcal{A})~\textmd{and}~\varrho(\mathcal{A})<P_s^t(\mathcal{A})$, we have
$\varrho(\mathcal{A})<\min\{P_t^s(\mathcal{A}),P_s^t(\mathcal{A})\}.$
Furthermore,
\begin{eqnarray}\label{th3-equ1}
\varrho(\mathcal{A})\leq \max\limits_{i,j\in N, j\neq i}\min\{P_i^j(\mathcal{A}),P_j^i(\mathcal{A})\}.
\end{eqnarray}

If
$\varrho(\mathcal{A})\in \Big(\tilde{\Psi}_{t,s}(\mathcal{A})\bigcap \mathcal{K}_t(\mathcal{A})\Big),$
i.e.,
$\varrho(\mathcal{A})\leq R_t(\mathcal{A})$
and
\begin{eqnarray}\label{th3-equ2}
\big(\varrho(\mathcal{A})-P_t^s(\mathcal{A})\big)\big(\varrho(\mathcal{A})-P_s^t(\mathcal{A})\big)
\leq \big(R_t(\mathcal{A})-P_t^s(\mathcal{A})\big)\big(R_s(\mathcal{A})-P_s^t(\mathcal{A})\big),
\end{eqnarray}
then solving $\varrho(\mathcal{A})$ in (\ref{th3-equ2}) gives
\begin{eqnarray*}
\varrho(\mathcal{A})\leq\frac{1}{2}\left\{
P_t^s(\mathcal{A})+P_s^t(\mathcal{A})+
\sqrt{\big(P_t^s(\mathcal{A})-P_s^t(\mathcal{A})\big)^2+4\big(R_t(\mathcal{A})-P_t^s(\mathcal{A})\big)\big(R_s(\mathcal{A})-P_s^t(\mathcal{A})\big)}
\right\}=\Delta_{t,s}(\mathcal{A}),
\end{eqnarray*}
and furthrermore
\begin{eqnarray}\label{th3-equ3}
\varrho(\mathcal{A})\leq\min\left\{R_t(\mathcal{A}),\Delta_{t,s}(\mathcal{A})\right\}
\leq\max\limits_{i,j\in N,j\neq i}\min\left\{R_i(\mathcal{A}),\Delta_{i,j}(\mathcal{A})\right\}.
\end{eqnarray}
The conclusion follows from (\ref{th3-equ1}) and (\ref{th3-equ3}).
\hfill$\blacksquare$

By Theorem \ref{th2}, Theorem 4.6 and Corollary 4.2 in \cite{wg}, the following comparison theorem can be derived easily.

\begin{theorem}\label{th4}
Let $\mathcal{A}=(a_{i_{1}\cdots i_{m}})\in{\mathbb{R}}^{[m,n]}$ be a weakly symmetric nonnegative tensor. Then
the upper bound in Theorem \ref{th3} is sharper than those in Theorem 4.6 of \cite{wg} and Corollary 4.5 of \cite{sys},
that is,
\begin{eqnarray*}
\varrho(\mathcal{A})&\leq&\Omega_{max}\\
&\leq&\max\limits_{i,j\in N,i\neq j}\left\{
\frac{1}{2}\Big(R_i(\mathcal{A})-a_{ij\cdots j}+P_j^i(\mathcal{A})+\Lambda^{\frac{1}{2}}(\mathcal{A})\Big),
R_i(\mathcal{A})-a_{ij\cdots j},P_j^i(\mathcal{A})\right\}\\
&\leq&\max\limits_{i\in N}R_i(\mathcal{A}),
\end{eqnarray*}
where
\begin{eqnarray*}
\Lambda_{i,j}(\mathcal{A})=(R_i(\mathcal{A})-a_{ij\cdots j}-P_j^i(\mathcal{A}))^2+4a_{ij\cdots j}(R_j(\mathcal{A})-P_j^i(\mathcal{A})).
\end{eqnarray*}
\end{theorem}

Finally, we show that the upper bound in Theorem \ref{th3} is sharper than those in \cite{wg,sys,liwen,hejunjcaa,hejunspringerplus,lql,hjaml} in some cases by the following two examples.
\begin{example}\label{example1}\emph{
Let $\mathcal{A}=(a_{ijkl})\in {\mathbb{R}}^{[4,2]}$ be a symmetric tensor defined by
$$a_{1111}=\frac{1}{2},~a_{2222}=3,~a_{ijkl}=\frac{1}{3}~ \textmd{elsewhere}.$$
By Corollary 4.5 of \cite{sys}, we have $$\varrho(\mathcal{A})\leq 5.3333.$$
By Theorem 2.7 of \cite{hjaml}, we have $$\varrho(\mathcal{A})\leq 5.2846.$$
By Theorem 3.3 of \cite{liwen}, we have $$\varrho(\mathcal{A})\leq 5.1935.$$
By Theorem 4.5, Theorem 4.6 and Theorem 4.7 of \cite{wg}, we all have $$\varrho(\mathcal{A})\leq 5.1822.$$
By Theorem 3.5 of \cite{hejunjcaa} and Theorem 6 of \cite{hejunspringerplus}, we both have $$\varrho(\mathcal{A})\leq 5.1667.$$
By Theorem 2.9 of \cite{lql}, we have $$\varrho(\mathcal{A})\leq 4.5147.$$
By Theorem \ref{th3}, we obtain $$\varrho(\mathcal{A})\leq 4.3971.$$}
\end{example}

\begin{example}\label{example2}\emph{
Let $\mathcal{A}=(a_{ijk})\in {\mathbb{R}}^{[3,3]}$ with entries defined as follows:
$$\mathcal{A}(:,:,1)=\left(
\begin{array}{ccc}
0&3&3\\
2.5&1&1\\
3&1&0\\
\end{array}\right),
~\mathcal{A}(:,:,2)=\left(
\begin{array}{ccc}
2&0.5&1\\
0&2&0\\
1&0.5&0\\
\end{array}\right),
~\mathcal{A}(:,:,3)=\left(
\begin{array}{ccc}
3&1&1\\
1&1&0\\
2&0&1\\
\end{array}\right).
$$
It is not difficult to verify that $\mathcal{A}$ is a weakly symmetric nonnegative tensor.
By Corollary 4.5 of \cite{sys} and Theorem 3.3 of \cite{liwen}, we both have $$\varrho(\mathcal{A})\leq 14.5000.$$
By Theorem 3.5 of \cite{hejunjcaa}, we have $$\varrho(\mathcal{A})\leq 14.2650.$$
By Theorem 4.6 of \cite{wg}, we have $$\varrho(\mathcal{A})\leq 14.2446.$$
By Theorem 4.5 of \cite{wg}, we have $$\varrho(\mathcal{A})\leq 14.1027.$$
By Theorem 6 of \cite{hejunspringerplus}, we have $$\varrho(\mathcal{A})\leq 14.0737.$$
By Theorem 4.7 of \cite{wg}, we have $$\varrho(\mathcal{A})\leq 13.2460.$$
By Theorem 2.9 of \cite{lql}, we have $$\varrho(\mathcal{A})\leq 13.2087.$$
By Theorem \ref{th3}, we obtain $$\varrho(\mathcal{A})\leq 11.7268.$$
}\end{example}

\begin{remark}\emph{
It is easy to see that in some cases the upper bound in Theorem \ref{th3} is sharper than those in [4-10] from Example \ref{example1} and Example \ref{example2}.
}\end{remark}

\section{Conclusions}

In this paper, we establish a new $Z$-eigenvalue localization set $\Omega(\mathcal{A})$ and prove that this set is tighter than those in \cite{wg}.
As an application, we obtain a new upper bound $\Omega_{max}$ for the $Z$-spectral radius of weakly symmetric nonnegative tensors,
and show that this bound is sharper than those in \cite{wg,sys,liwen,hejunjcaa,hejunspringerplus,lql,hjaml} in some cases by two numerical examples.

\section*{Competing interests}
The authors declare that they have no competing interests.
\section*{Authors' contributions}
All authors contributed equally to this work. All authors read and approved the final manuscript.
\section*{Acknowledgments}
This work is supported by the National Natural Science Foundations of China (Grant Nos.11361074,11501141),
Foundation of Guizhou Science and Technology Department (Grant No.[2015]2073)
and Natural Science Programs of Education Department of Guizhou Province (Grant No.[2016]066).



\bibliographystyle{elsarticle-num}

\end{document}